\documentclass[12pt]{amsart}
\usepackage[colorlinks=true,citecolor=blue,linkcolor=blue]{hyperref}
\usepackage{amssymb,amscd,amsmath,graphicx,mathtools}
\usepackage{enumerate}
\usepackage{fullpage}
\usepackage{comment}

\usepackage[all]{xy}
\usepackage{cleveref}

\numberwithin{equation}{section}
\newtheorem{theorem}{Theorem}[section]
\newtheorem{lemma}[theorem]{Lemma}

\newtheorem{corollary}[theorem]{Corollary}

\theoremstyle{definition}
\newtheorem{definition}[theorem]{Definition}

\newtheorem{question}[theorem]{Question}

\theoremstyle{definition}
\newtheorem{remark}[theorem]{Remark}
\newtheorem{example}[theorem]{Example}

\DeclareMathOperator{\Ass}{Ass}

\DeclareMathOperator{\NP}{NP}
\DeclareMathOperator{\SP}{SP}

\DeclareMathOperator{\conv}{convexhull}

\newcommand{\PP}{{\mathbb P}}
\newcommand{\ZZ}{{\mathbb Z}}
\newcommand{\NN}{{\mathbb N}}
\newcommand{\QQ}{{\mathbb Q}}
\newcommand{\RR}{{\mathbb R}}

\def\B{{\mathcal B}}
\def\C{{\mathcal C}}

\def\R{{\mathcal R}}

\def\R{{\mathcal R}}

\def\h{\widetilde{H}}

\def\aa{{\mathfrak a}}

\def\pp{{\mathfrak p}}

\def\bb{{\mathfrak b}}

\def\a{{\bf a}}
\def\b{{\bf b}}

\def\h{{\bf h}}

\def\v{{\bf v}}
\def\x{{\bf x}}
\def\y{{\bf y}}

\newcommand{\kk}{\Bbbk}

\def\rhat{\widehat{\rho}}

\def\vhat{\widehat{v}}

\def\1{{\bf 1}}
\def\0{{\bf 0}}

\begin{document}
\title[Resurgence and convex body of graded families of ideals]{Resurgence number and convex body associated to pairs of graded families of ideals}


%    author one information
\author{T\`ai Huy H\`a}
\address{Tulane University, Department of Mathematics, 6823 St. Charles Ave., New Orleans, LA 70118, USA}
\email{tha@tulane.edu}
%\urladdr{???}

\author{A.V. Jayanthan}
\address{Department of Mathematics, I.I.T. Madras, Chennai - 600036, INDIA}
\email{jayanav@iitm.ac.in}

\author{Arvind Kumar}
\address{Department of Mathematical Sciences, New Mexico State University, 1305 Frenger St, Las Cruces, NM 88001, USA}
\email{arvkumar@nmsu.edu}

\author{Th\'ai Th\`anh Nguy$\tilde{\text{\^E}}$n}
\address{Department of Mathematics, University of Dayton, 300 College Park
Dayton, OH 45469, USA and \\
	University of Education, Hue University, 34 Le Loi St., Hue, Viet Nam}
\email{tnguyen11@tulane.edu ntthai.dhsp@hueuni.edu.vn}
\thanks{}

%    \subjclass is required.

\date{}

\dedicatory{}

%    "Communicated by" -- provide editor's name; required.
\commby{}

%    Abstract is required.
\begin{abstract}
We discuss how to understand the asymptotic resurgence number of a pair of graded families of ideals from combinatorial data of their associated convex bodies. When the families consist of monomial ideals, the convex bodies being considered are the Newton-Okounkov bodies of the families. When ideals in the second family are classical invariant ideals, for instance, determinantal ideals or ideals of Pfaffians, these convex bodies are constructed from the associated Rees packages.
\end{abstract}

\maketitle

\section{Introduction} \label{sec.intro}

Let $R$ be a Noetherian commutative ring, and let $\aa_\bullet = \{\aa_i\}_{i \ge 1}$ and $\bb_\bullet = \{\bb_i\}_{i \ge 1}$ be graded families of ideals in $R$. The \emph{resurgence} and \emph{asymptotic resurgence} numbers, $\rho(\aa_\bullet,\bb_\bullet)$ and $\rhat(\aa_\bullet,\bb_\bullet)$, are measures for the non-containment between elements in the families $\aa_\bullet$ and $\bb_\bullet$; see \cite{HKNN23}. These invariant are generalizations of the resurgence and asymptotic resurgence numbers of an ideal, that have been much investigated in the \emph{Ideal Containment Problem} and have shown many exciting applications (cf. \cite{BGHN22a, BGHN22b, BN24a, BN24b, BH, CHHVT, DFMS19, DD2020, DNS23, DHNSTT, GHM2020, GSV22, GHV, HKZ, JKM2021, Ng23a, Ng23b, Villarreal22} and references therein). 
Particularly,
\begin{align*}
	\rho(\aa_\bullet, \bb_\bullet) & = \sup\left\{ \dfrac{s}{r} ~\Big|~ s, r \ge 1, \aa_s \not\subseteq \bb_r \right\}, \text{ and } \\
	\rhat(\aa_\bullet,\bb_\bullet) & = \sup\left\{ \dfrac{s}{r} ~\Big|~ s, r \ge 1, \aa_{st} \not\subseteq \bb_{rt} \text{ for  } t \gg 0 \right\}.
\end{align*}
Resurgence and asymptotic resurgence a priori are very difficult to compute, even in the case for one ideal.

Our goal is to provide a combinatorial understanding of the asymptotic resurgence of a pair of graded families of ideals through associated convex bodies. Our first result addresses the case when the families $\aa_\bullet$ and $\bb_\bullet$ consist of monomial ideals. In this case, the convex bodies of interest are the \emph{Newton-Okounkov bodies} $\Delta(\aa_\bullet)$ and $\Delta(\bb_\bullet)$ associated to $\aa_\bullet$ and $\bb_\bullet$; see Section \ref{sec.NO} for precise definitions.
In general, the relationship between algebraic invariant and properties of a graded family of ideals and combinatorial data from its associated Newton-Okounkov body has been much studied (cf. \cite{Cut2013, Cut2014,HN23, KK2012, KK2014, LM2009} and references thereafter). Our work adds another layer to this extensive research program.

Leaving out the detailed definitions and terminology until later, our first main result is stated as follows. (Note that $\Delta(\bb_\bullet)^o$ is the \emph{polar set} of $\Delta(\bb_\bullet)$, as defined in Definition \ref{def.polar}.)

\medskip

\noindent \textbf{Theorem \ref{thm.rhohatViaNO}.} Let $\aa_\bullet$ and $\bb_\bullet$  be graded families of monomial ideals in $R=\kk[x_1,\ldots,x_n]$.
\begin{enumerate}
	\item If $\R(\bb_\bullet)$ is Noetherian, then
	$$\rhat(\aa_\bullet, \overline{\bb_\bullet}) =\sup \{ \lambda > 0 ~\big|~ \lambda \cdot \Delta(\aa_\bullet) \not\subseteq \Delta(\bb_\bullet) \}.$$
	\item If $\{ \lambda > 0 ~\big|~ \lambda \cdot \Delta(\aa_\bullet) \not\subseteq \Delta(\bb_\bullet)\} \not= \emptyset$, then
	$$\sup \{ \lambda > 0 ~\big|~ \lambda \cdot \Delta(\aa_\bullet) \not\subseteq \Delta(\bb_\bullet) \} = \dfrac{1}{\inf\{ \langle \mathbf{a}, \mathbf{b} \rangle \mid \mathbf{a}\in \Delta(\aa_\bullet) ,  \mathbf{b}\in \Delta(\bb_\bullet)^o\}}. $$
\end{enumerate}

\medskip

The statement of Theorem \ref{thm.rhohatViaNO} is inspired by an attempt to understand the main theorem of \cite{Villarreal22}. More precisely, for an ideal $I \subseteq R$, write $I^{(\bullet)} = \{I^{(i)}\}_{i \ge 1}$ and $\overline{I^\bullet} = \{\overline{I^i}\}_{i \ge 1}$ for the graded families of symbolic and integral closure of ordinary powers of $I$. \cite[Theorem 3.7]{Villarreal22} shows that if $I$ is a \emph{squarefree monomial} ideal in a polynomial ring $R = \kk[x_1, \dots, x_n]$, then
$$\dfrac{1}{\rho(I^{(\bullet)}, \overline{I^\bullet})} = \min\{ \langle u,v \rangle ~\big|~ u \in \SP(I), v \in \SP(I^\vee)\},$$
where $I^\vee$ is the Alexander dual of $I$, and $\SP(I)$ represents the \emph{symbolic polyhedron} of $I$. We will see in Corollary \ref{cor.rhatV} that this formula is equivalent to
\begin{align}
	\rhat(I^{(\bullet)}, \overline{I^\bullet}) = \sup \{\lambda > 0 ~\big|~ \lambda \cdot \SP(I) \not\subseteq \NP(I)\}; \label{eq.V}
\end{align}
here, $\NP(I)$ denotes the Newton polyhedron of $I$.
In many important situations, by \cite[Corollary 3.9]{HKNN23}, we have $\rhat(\aa_\bullet, \bb_\bullet) = \rhat(\aa_\bullet, \overline{\bb_\bullet}) = \rho(\aa_\bullet, \overline{\bb_\bullet})$; particularly, in a polynomial ring, $\rhat(I^{(\bullet)}, I^\bullet) = \rhat(I^{(\bullet)}, \overline{I^\bullet}) = \rho(I^{(\bullet)}, \overline{I^\bullet})$ --- see also \cite{DFMS19}.

\medskip

The problem is much harder for graded families of arbitrary ideals. To get similar statements for more general ideals, we shall focus on the class of \emph{Nagata rings}; those are Noetherian \emph{domains} $R$ satisfying the property that every finitely generated domain over $R$ has a module-finite integral closure in any finite extension of its quotient field, i.e., $R$ is \emph{universally Japanese}. The restriction to Nagata rings is to avoid strange behaviors, for instance, where $\R(\overline{\bb})$ is not a finitely generated $\R(\bb)$-module, as exhibited in Nagata's example (cf. \cite[Example 3.5]{HKNN23}).

Let $\bb$ be an ideal in a Nagata ring $R$ that has a \emph{Rees package} $\R = (\B,\v,\Gamma)$. The notion and existence of Rees package was introduced and studied recently in \cite{BDHM}. For a graded family $\aa_\bullet = \{\aa_i\}_{i \ge 1}$ of ideals in $R$, there associates a closed set $\Gamma_\R(\aa_\bullet)$; see Section \ref{sec.Rees} for the precise definitions. A graded family $\bb_\bullet = \{\bb_i\}_{i \ge 1}$ is called \emph{$\bb$-equivalent} if there exists an integer $k \in \NN$ such that for all $i \ge 1$,
$$\bb_{i+k} \subseteq \bb^i \subseteq \bb_i.$$
Our next result is stated as follows.

\medskip

\noindent \textbf{Theorems \ref{thm.rhat_b} and \ref{thm.rhat_b-equivalent}.} Let $\bb \subseteq R$ be an ideal in a Nagata ring $R$ with Rees package $\R = (\B,\v,\Gamma)$, in which $\v$ is a vector of valuations. Let $\aa_\bullet$ and $\bb_\bullet$ be graded families of nonzero ideals in $R$ such that either
\begin{enumerate}
	\item $\bb_\bullet = \bb^\bullet$ is the family of ordinary powers of $\bb$; or
	\item $\aa_\bullet$ is a filtration and $\bb_\bullet$ is $\bb$-equivalent.
\end{enumerate}
Then,
$$\rhat(\aa_\bullet, \overline{\bb_\bullet}) = \sup \{\lambda > 0 ~\big|~ \lambda \cdot \Gamma_\R(\aa_\bullet) \not\subseteq \Gamma\}.$$

\medskip

Results in \cite{BDHM} show that the following classes of ideals possess Rees packages:
\begin{enumerate}
	\item[(i)] monomial ideals in affine semigroup rings;
	\item[(ii)] sums of products of determinantal ideals of generic matrices;
	\item[(iii)] sums of products of determinantal ideals of symmetric matrices;
	\item[(iv)] sums of products of ideals of Pfaffians of skew-symmetric matrices;
	\item[(v)] products of determinantal ideals of Hankel matrices.
\end{enumerate}
Thus, Theorems \ref{thm.rhat_b} and \ref{thm.rhat_b-equivalent} are applicable when $\bb \subseteq R$ is an ideal belonging to the classes (i) -- (v). In these cases, $\Gamma_\R(\aa_\bullet)$ is a closed convex set.

\medskip

Our study in Sections \ref{sec.NO} and \ref{sec.Rees} assumes the graded family $\bb_\bullet$ to be Noetherian; that is, the Rees algebra $\R(\bb_\bullet)$ is a Noetherian ring. In general, $\aa_\bullet$ and $\bb_\bullet$ are not necessarily Noetherian. For a non-Noetherian graded family $\aa_\bullet$, one can approximate it by a collection of Noetherian families, the \emph{truncations} of $\aa_\bullet$; see Section \ref{sec.Truncate} for precise definitions. Our last main result shows that resurgence and asymptotic resurgence numbers, $\rho(\aa_\bullet, \bb_\bullet)$ and $\rhat(\aa_\bullet,\bb_\bullet)$, can be computed via those of truncations of $\aa_\bullet$. On the other hand, examples exist to exhibit that truncations of $\bb_\bullet$ will not give the same conclusion; see Example \ref{ex.rhoTruncation}.

\medskip

\noindent \textbf{Theorem \ref{thm.resurgencetruncation1}.}  Let $\aa_\bullet$ and $\bb_\bullet$ be graded families of ideals in a Noetherian commutative ring $R$. For $n \in \NN$, let $\aa_{n,\bullet}$ denote the $n$-th truncation graded family of $\aa_\bullet$. We have
$$\rho(\aa_\bullet,\bb_\bullet) =\sup\limits_{n \ge 1 } \rho(\aa_{n,\bullet},\bb_\bullet)=\lim_{n \to \infty} \rho(\aa_{n,\bullet},\bb_\bullet) \text{ and }
\rhat(\aa_\bullet,\bb_\bullet) =\sup\limits_{n \ge 1 } \rhat(\aa_{n,\bullet},\bb_\bullet)=\lim_{n \to \infty} \rhat(\aa_{n,\bullet},\bb_\bullet).$$

\medskip

%%%%%%%%%%%%%%%%%%%%%%%%%%%%%%%%

\section{Asymptotic resurgence and Newton-Okounkov body} \label{sec.NO}

This section is devoted to demonstrating a connection between the asymptotic resurgence number $\rhat(\aa_\bullet, \overline{\bb_\bullet})$ of a pair of graded families of ideals and Newton-Okounkov bodies associated with these families.

We collect some important notations and terminology. If $R = \kk[x_1, \dots, x_n]$ is a polynomial ring over a field $\kk$ and $\a = (a_1, \dots, a_n) \in \ZZ_{\ge 0}^n$, then we shall denote by $x^\a$ the monomial $x_1^{a_1} \dots x_n^{a_n}$ in $R$. The Newton polyhedron of a monomial ideal is a well-known notion. The construction of a symbolic polyhedron for a monomial ideal is from \cite{CEHH2017}.

\begin{definition}
	\label{def.NP}
	Let $R = \kk[x_1, \dots, x_n]$ and let $I \subseteq R$ be a monomial ideal.
	\begin{enumerate}
		\item The \emph{Newton polyhedron} of $I$ is defined to be
	$$\NP(I) = \conv \langle \{ \a \in \ZZ_{\ge 0}^n ~\big|~ x^\a \in I \}\rangle \subseteq \RR^n_{\ge 0}.$$
		\item Suppose, in addition, that $I$ is squarefree. The \emph{symbolic polyhedron} of $I$ is given by
		$$\SP(I) = \bigcap_{\pp \in \Ass(I)} \NP(\pp) = \bigcap_{\pp \in \Ass(I)} \left\{(a_1, \dots, a_n) \in \RR_{\ge 0}^n ~\big|~ \sum_{x_i \in \pp} a_i \ge 1\right\} \subseteq \RR_{\ge 0}^n.$$
	\end{enumerate}
\end{definition}

We introduce a notion of polar sets that contrasts with the commonly known concept in analysis. For clarity and simplicity, we will still refer to them as polar sets throughout this paper, avoiding the need for new terminology.

\begin{definition}
	\label{def.polar}
For a convex set $\C \subseteq \RR^n$ not containing the origin, the \emph{polar set} of $\C$ is defined by
	$$\C^o = \{\a \in \RR^n ~\big|~ \langle \a, \b \rangle \ge 1 \text{ for all } \b \in \C\},$$
	where $\langle \bullet, \bullet \rangle$ denotes the usual dot product in $\RR^n$.
\end{definition}
	
\begin{remark}
	\label{rmk.dual}
	From the discussion in \cite{DNS2, GRV09, GSV22}, it can be observed that, for a squarefree monomial ideal $I \subseteq \kk[x_1, \dots, x_n]$,
	$$\SP(I)^o = \NP(I^\vee) \text{ and } \NP(I) = \SP(I^\vee)^o,$$
	where $I^\vee = \langle\prod_{x_i \in \pp} x_i ~\big|~ \pp \in \Ass(I) \rangle$ denotes the \emph{Alexander dual} of $I$.
\end{remark}

The ``Bipolar Theorem'' is a well-known result in analysis (cf. \cite{Rock}). Since our notion of polar sets differs from the analytical definition, we will provide the proof of a similar bipolar statement in a simple case relevant to our study.

\begin{lemma}
	\label{lem.bipolar}
	Let $\C \subseteq \RR^n_{\ge 0}$ be a closed convex set that does not contain the origin and absorbs $\RR^n_{\ge 0}$, i.e., $\C + \RR^n_{\ge 0} \subseteq \C$. Then,
	$$\C^{oo} = C.$$
\end{lemma}

\begin{proof}
	It is easy to see from the definition that $\C \subseteq \C^{oo}$. We shall establish the reverse containment.
	Suppose, by contradiction, that there exists $\mathbf{b} \in \C^{oo} \setminus \C$.
	
	It is a standard result from convex geometry that there exists a hyperplane $H$ which \emph{strongly} separates $\mathbf{b}$ and $\C$, whose equation is of the form $\langle \mathbf{c}, \x \rangle = 1$, where $\mathbf{c} = (c_1, \dots, c_n)$.  Since $\C \subseteq \RR^n_{\ge 0}$, we may assume that $H$ contains points in $\RR^n_{\ge 0}$ which are not the origin. This implies that there exists $i$ such that $c_i > 0$.
	
	Since $\C$ absorbs $\RR^n_{\ge 0}$, by increasing the $i$-th coordinate of the points in $\C$, it follows that $\C$ lies in the positive half-space determined by $H$. That is, $\langle \mathbf{c},\y\rangle \ge 1$ for all $\y \in \C$. This implies that $\mathbf{c} \in \C^o$.
	On the other hand, since $\b$ lies on the other side of $H$, we must have $\langle \mathbf{c}, \mathbf{b}\rangle < 1$. As a consequence, $\b \not\in (\C^o)^o = \C^{oo}$, which is a contradiction.
\end{proof}

A \emph{graded family} of ideals in a commutative ring $R$ is a collection of ideals $\{\aa_i\}_{i \ge 1}$ satisfying the condition that $\aa_p \cdot \aa_q \subseteq \aa_{p+q}$ for all $p, q \ge 1$. A graded family $\{\aa_i\}_{i \ge 1}$ is called a \emph{filtration} if $\aa_i \supseteq \aa_{i+1}$ for all $i \ge 1$. The \emph{Rees algebra} associated to a graded family $\aa_\bullet = \{\aa_i\}_{i \ge 1}$ of ideals in $R$ is defined to be
$$\R(\aa_\bullet) = \bigoplus_{i \ge 0} \aa_i t^i \subseteq R[t],$$
where, by convention, $\aa_0 = R$. The graded family $\aa_\bullet$ is called \emph{Noetherian} if the Rees algebra $\R(\aa_\bullet)$ is a Noetherian ring.

The \emph{Newton-Okounkov body} associated with a graded family of ideals, in general, is defined based on the existence of \emph{good valuations}; see \cite{Cut2013, Cut2014, KK2012, KK2014}. This construction becomes more transparent for graded families of monomial ideals; see, for instance, \cite{HN23, KK2014}.

\begin{definition}
	\label{def.NObody}
	Let $R = \kk[x_1, \dots, x_n]$ and let $\aa_\bullet = \{\aa_i\}_{i \ge 1}$ be a graded family of \emph{monomial} ideals in $R$. The \emph{Newton-Okounkov body} associated to $\aa_\bullet$ is
	$$\Delta(\aa_\bullet) = \overline{\bigcup_{k \ge 1} \dfrac{1}{k} \NP(\aa_k)} \subseteq \RR_{\ge 0}^n.$$
\end{definition}

\begin{remark}
Observe that if $\aa_\bullet = I^{\bullet}$, for a monomial ideal $I$, then $\Delta(\aa_\bullet) = \NP(I)$, and if $\aa_\bullet = I^{(\bullet)}$, for a squarefree monomial ideal $I$, then $\Delta(\aa_\bullet) = \SP(I)$.

By construction, $\Delta(\aa_\bullet)$ is a closed convex set that absorbs $\RR^n_{\ge 0}$. It was also pointed out in \cite[Example 3.2]{HN23} that for any nonempty closed convex set $P \subseteq \RR_{\ge 0}^n$ absorbing $\RR_{\ge 0}^n$, there is a graded family of monomial ideal $\aa_\bullet$ such that $\Delta(\aa_\bullet) = P$.
\end{remark}

We are now ready to present our first main result.

\begin{theorem}
\label{thm.rhohatViaNO}
Let $\aa_\bullet$ and $\bb_\bullet$  be graded families of monomial ideals in $R=\kk[x_1,\ldots,x_n]$.
\begin{enumerate}
	\item If $\R(\bb_\bullet)$ is Noetherian, then
$$\rhat(\aa_\bullet, \overline{\bb_\bullet}) =\sup \{ \lambda > 0 ~\big|~ \lambda \cdot \Delta(\aa_\bullet) \not\subseteq \Delta(\bb_\bullet) \}.$$
	\item If $\{ \lambda > 0 ~\big|~ \lambda \cdot \Delta(\aa_\bullet) \not\subseteq \Delta(\bb_\bullet)\} \not= \emptyset$, then
	$$\sup \{ \lambda > 0 ~\big|~ \lambda \cdot \Delta(\aa_\bullet) \not\subseteq \Delta(\bb_\bullet) \} = \dfrac{1}{\inf\{ \langle \mathbf{a}, \mathbf{b} \rangle \mid \mathbf{a}\in \Delta(\aa_\bullet) ,  \mathbf{b}\in \Delta(\bb_\bullet)^o\}}. $$
\end{enumerate}
\end{theorem}

\begin{proof} (1) We will first consider the trivial case, when $\rhat(\aa_\bullet, \overline{\bb_\bullet}) = -\infty$. By definition, this happens if only if for any integers $s, r \ge 1$ and any $M \in \NN$, there exists $\theta \ge M$ such that $\aa_{s\theta} \subseteq \overline{\bb_{r\theta}}$. Particularly, $\overline{\aa_{s\theta}} \subseteq \overline{\bb_{r\theta}}$, and so
$$\frac{s}{r} \cdot \frac{1}{s\theta}\NP(\aa_{s\theta}) \subseteq \frac{1}{r\theta}\NP(\bb_{r\theta}) \subseteq \Delta(\bb_\bullet).$$
 Thus, for any $\lambda = \frac{p}{q} \in \QQ_{> 0}$ and any $\ell \in \NN$, by taking $s = p\ell$ and $r = q\ell$, we obtain
$$\lambda \cdot \dfrac{1}{\ell}\NP(\aa_\ell) = \dfrac{p\ell}{q\ell} \cdot \dfrac{1}{\ell}\NP(\aa_\ell) \subseteq \dfrac{p\ell}{q\ell} \cdot \dfrac{1}{p\ell \theta} \NP(\aa_{p\ell \theta}) = \dfrac{s}{r} \cdot \dfrac{1}{s\theta} \NP(\aa_{s\theta}) \subseteq \Delta(\bb_\bullet).$$
Therefore, $\lambda \cdot \Delta(\aa_\bullet) \subseteq \Delta(\bb_\bullet)$ for any $\lambda \in \QQ_{> 0}$. Hence, we get $\lambda \cdot \Delta(\aa_\bullet) \subseteq \Delta(\bb_\bullet)$ for any $\lambda > 0$, i.e.,
$$\sup \{\lambda > 0 ~\big|~ \lambda \cdot \Delta(\aa_\bullet) \not\subseteq \Delta(\bb_\bullet)\} = -\infty.$$

Conversely, suppose that for any $\lambda > 0$, we have $\lambda \cdot \Delta(\aa_\bullet) \subseteq \Delta(\bb_\bullet)$. Since $\R(\bb_\bullet)$ is Noetherian, it follows from \cite[Theorem 3.4]{HN23} that there exists $k \ge 1$ such that $\overline{\bb_k^t} = \overline{\bb_{kt}}$ and $\Delta(\bb_\bullet)=\frac{1}{k} \NP(\bb_k)= \frac{1}{kt} \NP(\bb_{kt})$ for all $t \in \NN$.
As a consequence, for any integers $s,r \in 1$, we have
$$\frac{s}{r} \cdot \frac{1}{skt}\NP(\aa_{skt}) \subseteq \Delta(\bb_\bullet) = \frac{1}{rkt}\NP(\bb_{rkt}) \text{ for all } t \in \NN.$$
Particularly, $\NP(\aa_{skt}) \subseteq \NP(\bb_{rkt})$. It then follows that $\aa_{skt} \subseteq \overline{\aa_{skt}} \subseteq \overline{\bb_{rkt}}$. Therefore, for any $M \in \NN$, by choosing $t > \frac{M}{k}$ and setting $\theta = kt > M$, we have $\aa_{s\theta} \subseteq \overline{\bb_{r\theta}}$. Hence,
$$\rhat(\aa_\bullet, \overline{\bb_\bullet}) = -\infty.$$

Consider the case that $\rhat(\aa_\bullet, \overline{\bb_\bullet}) \not= -\infty$. Let $s,r \ge 1$ be such that $\aa_{st} \not\subseteq \overline{\bb_{rt}}$ for $t\gg 0$; the existence of $s$ and $r$ are due to the assumption on $\rhat(\aa_\bullet, \overline{\bb_\bullet})$.
Particularly, for $t \gg 0$, there exists a monomial $x_1^{u_{t,1}} \cdots x_n^{u_{t,n}} \in \aa_{st} \setminus \overline{\bb_{rt}}$. This implies that
$u_t =(u_{t,1},\ldots,u_{t,n}) \in \NP(\aa_{st}) \setminus \NP(\bb_{rt}).$

By replacing $t$ with $kt$, it follows that for $t \gg 0$,
$$u_{kt} =(u_{kt,1},\ldots,u_{kt,n}) \in \NP(\aa_{skt}) \setminus \NP(\bb_{rkt}).$$
We now have
$$\frac{u_{kt}}{skt} \in \frac{1}{skt} \NP(\aa_{skt}) \subseteq \Delta(\aa_\bullet) \text{ and } \frac{s}{r} \cdot \frac{u_{kt}}{skt}=\frac{u_{kt}}{rkt}  \not\in\frac{1}{rkt} \NP(\bb_{rkt})=\Delta(\bb_\bullet).$$
Therefore, $\frac{s}{r} \cdot \Delta(\aa_\bullet) \not\subseteq \Delta(\bb_\bullet)$ and, hence,
$$\rhat(\aa_\bullet, \overline{\bb_\bullet})\le \sup \{ \lambda > 0 ~\big|~  \lambda \cdot \Delta(\aa_\bullet) \not\subseteq \Delta(\bb_\bullet) \} .$$

Conversely, let $s,r \ge 1$ be such that $\frac{s}{r} \cdot \Delta(\aa_\bullet) \not\subseteq \Delta(\bb_\bullet).$ Then, there exists $q \in \NN$ such that
$$\frac{s}{r} \cdot \frac{1}{q} \NP(\aa_q) \not\subseteq \Delta(\bb_\bullet) =\frac{1}{kt} \NP(\bb_{kt}) \text{ for all } t \in \NN.$$
This implies that there is a rational vector $v=(v_1,\ldots,v_n) \in \NP(\aa_q)$ such that $\frac{s}{r} \cdot \frac{1}{q} v \not\in \frac{1}{kt} \NP(\bb_{kt})$ for all $t$. Particularly,
$$\frac{s}{r} \cdot \frac{1}{q} v \not\in \frac{1}{kq}\NP(\bb_{kq}).$$
By replacing $v$ and $q$ with $vL$ and $qL$, respectively, where $L$ is the least common multiple of the denominators of $v_1, \dots, v_n$, we may assume that $v$ is an integral vector.

Since $v \in \NP(\aa_q)$, $x_1^{v_1} \ldots x_n^{v_n} \in \overline{\aa_q}$, which further implies that $x_1^{pv_1} \ldots x_n^{pv_n} \in \aa_q^p \subseteq \aa_{pq}$ for some $p \in \NN$.
We then have
$$\left(x_1^{pv_1} \ldots x_n^{pv_n}\right)^{skt} \in \aa_{pq}^{skt} \subseteq \aa_{sktpq} \text{ and } sktp \cdot v \not\in rtp \cdot \NP(\bb_{kq}).$$
It follows that
$$(x^v)^{sktp} = \left(x_1^{pv_1} \ldots x_n^{pv_n}\right)^{stk} \in \aa_{sktpq} \setminus  \overline{\bb_{kq}^{rtp}}= \aa_{sktpq} \setminus \overline{\bb_{rtpkq}} .$$
Thus, $\aa_{sktpq} \not\subseteq \overline{\bb_{rktpq}}$ for all $t \in \NN$. Therefore, $\frac{s}{r} =\frac{skpq}{rkpq} \le \rhat(\aa_\bullet,\overline{\bb_\bullet}).$ Hence,
$$\sup \{ \lambda > 0 ~\big|~ \lambda \cdot \Delta(\aa_\bullet) \not\subseteq \Delta(\bb_\bullet) \}  \le \rhat(\aa_\bullet,\overline{\bb_\bullet}),$$
and the first assertion of the theorem is proved.

(2) To establish the second equality, we use the fact $\lambda \cdot \Delta(\aa_\bullet) \subseteq \gamma \cdot \Delta(\aa_\bullet)$ for any $0 < \gamma \le \lambda$ from \cite{DNS2}. It follows that
\begin{eqnarray*}
		\sup\{\lambda > 0 \mid \lambda \cdot \Delta(\aa_\bullet) \not\subseteq \Delta(\bb_\bullet)\}
		&=&\inf\{\lambda > 0 \mid \lambda \cdot \Delta(\aa_\bullet) \subseteq \Delta(\bb_\bullet)\}.
\end{eqnarray*}
We continue to show that
$$\inf\{\lambda > 0 \mid \lambda \cdot \Delta(\aa_\bullet) \subseteq \Delta(\bb_\bullet)\} = \inf\{\lambda > 0 \mid \langle \lambda \mathbf{a},\mathbf{b}\rangle \geq 1, \forall \mathbf{a}\in\Delta(\aa_\bullet), \mathbf{b}\in\Delta(\bb_\bullet)^o\}.$$
Indeed, if $\lambda \cdot \Delta(\aa_\bullet) \subseteq \Delta(\bb_\bullet)$, then for any $\mathbf{a} \in \Delta(\aa_\bullet)$ (whence $\lambda \mathbf{a} \in \Delta(\bb_\bullet)$) and $\mathbf{b} \in \Delta(\bb_\bullet)^o$, by definition, $\langle \lambda \mathbf{a}, \mathbf{b}\rangle \ge 1$.
On the other hand, if $\langle \lambda \mathbf{a}, \mathbf{b}\rangle \ge 1$ for all $\mathbf{a} \in \Delta(\aa_\bullet)$ and $\mathbf{b} \in \Delta(\bb_\bullet)^o$, then $\lambda \cdot \Delta(\aa_\bullet) \subseteq \Delta(\bb_\bullet)^{oo} = \Delta(\bb_\bullet)$, where the last equality follows from Lemma \ref{lem.bipolar}.

Finally, it follows that
\begin{eqnarray*}
		\sup\{\lambda > 0 \mid \lambda \cdot \Delta(\aa_\bullet) \not\subseteq \Delta(\bb_\bullet)\} &=& \inf\{\lambda > 0 \mid \langle \lambda \mathbf{a},\mathbf{b}\rangle \geq 1, \forall \mathbf{a}\in\Delta(\aa_\bullet), \mathbf{b}\in\Delta(\bb_\bullet)^o\}\\
		& = &\inf\left\{\lambda ~\Big|~ \dfrac{1}{\langle\mathbf{a},\mathbf{b}\rangle} \leq \lambda, \forall \mathbf{a}\in\Delta(\aa_\bullet), \mathbf{b}\in\Delta(\bb_\bullet)^o\right\}\\
		&=&\dfrac{1}{\inf\{ \langle \mathbf{a}, \mathbf{b} \rangle \mid \mathbf{a}\in \Delta(\aa_\bullet) ,  \mathbf{b}\in \Delta(\bb_\bullet)^o\}}.
	\end{eqnarray*}
The assertion is established.
\end{proof}

\begin{remark}
	In a different context, for graded families of polyhedra, a similar looking statement to that of Theorem \ref{thm.rhohatViaNO} is obtained in \cite{DNS2}.
\end{remark}

We illustrate Theorem \ref{thm.rhohatViaNO} with the following example. This example appeared in \cite{HN23} for a different purpose.

\begin{example}  Let $I =(x,y)^2 \cap (y,z)^3 \cap (x,z)^4 \subseteq \kk[x,y,z]$. Consider the graded families $\aa_\bullet = I^{(\bullet)}$ and $\bb_\bullet = I^\bullet$.  Then, it follows from \cite[Remark 2.14]{HN23} that $\Delta(\aa_\bullet)=\SP(I)$ and $\Delta(\bb_\bullet)=\NP(I).$ It is known from \cite[Example 5.5]{HN23} that  $$\SP(I) =\text{convexhull}\left\{ (4,3,0), (2,0,3), (0,2,4), \left( \frac{3}{2}, \frac{1}{2}, \frac{5}{2}\right)  \right\} +\mathbb{R}^3_{\ge 0}.$$ Note that $I= ( xyz^3,x^2z^3,x^2yz^2, y^2z^4,x^3y^2z,x^4y^3).$ Therefore, $$\NP(I)= \text{convexhull}\left\{ (1,1,3), (2,0,3), (2,1,2), (0,2,4), (3,2,1), (4,3,0)\right\}+\mathbb{R}^3_{\ge 0}.$$

Observe that the first three vertices of $\SP(I)$ belong to $\NP(I)$. Thus,
$$\sup\{ \lambda > 0 ~\big|~ \lambda \cdot \SP(I) \not\subseteq \NP(I)\} = \sup \left\{\lambda > 0 ~\Big|~ \lambda \cdot \left(\frac{3}{2}, \frac{1}{2}, \frac{5}{2}\right) \not\in \NP(I)\right\}.$$
This value is then determined by where the line connecting the origin and $\left(\frac{3}{2}, \frac{1}{2}, \frac{5}{2}\right)$ \emph{first} intersects $\NP(I)$. This point can be computed to be $\left(\frac{5}{3}, \frac{5}{9}, \frac{25}{9}\right)$. Thus,
$$\sup \left\{\lambda > 0 ~\Big|~ \lambda \cdot \left(\frac{3}{2}, \frac{1}{2}, \frac{5}{2}\right) \not\in \NP(I)\right\} = \frac{10}{9}.$$
Hence, due to Theorem \ref{thm.rhohatViaNO}, we get
$\rhat(\aa_\bullet, \overline{\bb_\bullet}) =\sup \{ \lambda > 0 ~\big|~ \lambda \cdot \Delta(\aa_\bullet) \not\subseteq \Delta(\bb_\bullet) \}=\frac{10}{9}.$

 It is worth noting that, in \cite[Theorem 3]{JZ21}, the resurgence number was also given as $\rho(\aa_\bullet, \bb_\bullet) = \frac{10}{9}$. Thus, there are equalities
 $$\rhat(\aa_\bullet, \overline{\bb_\bullet})= \rhat(\aa_\bullet, \bb_\bullet) = \rho(\aa_\bullet, \overline{\bb_\bullet}) = \rho(\aa_\bullet, \bb_\bullet) = \frac{10}{9}.$$
\end{example}

\begin{question}
	Does the statement of Theorem \ref{thm.rhohatViaNO}(1) still hold without the condition that $\R(\bb_\bullet)$ is Noetherian? We do not know of any example suggesting otherwise.
\end{question}

%%%%%%%%%%%

As a consequence of Theorem \ref{thm.rhohatViaNO}, we obtain the following generalization of \cite[Theorem 3.7]{Villarreal22}.

\begin{corollary} \label{cor:asym-res-IP}
    Let $\aa_\bullet$ be graded family  of monomial ideals in $R=\kk[x_1,\ldots,x_n]$ and $\bb$ be a squarefree monomial ideal in $R$. Then, $$\rhat(\aa_\bullet, \overline{\bb^\bullet}) =\sup \{ \lambda > 0 ~\big|~ \lambda \cdot \Delta(\aa_\bullet) \not\subseteq \NP(\bb) \} = \dfrac{1}{\inf\{ \langle \mathbf{a}, \mathbf{b} \rangle \mid \mathbf{a}\in \Delta(\aa_\bullet) ,  \mathbf{b}\in \SP(\bb^\vee)\}} . $$
\end{corollary}

\begin{proof}
By Remark \ref{rmk.dual} (see also \cite{DNS2,GRV09}), we have
$$\NP(\bb)= \SP(\bb^\vee)^o = \left\{ \a \in \RR^n_{\ge 0} ~\big|~ \langle \a, \b \rangle \ge 1 \text{ for all } \b \in \SP(\bb^\vee)\right\}.$$
Thus, applying Theorem \ref{thm.rhohatViaNO}, we get the desired result.
\end{proof}

%%%%%%%%

\cite[Theorem 3.7]{Villarreal22} is now recovered by letting $\aa_\bullet = \aa^{(\bullet)}$ for an ideal $\aa$, and noting that $\rhat(\aa^{(\bullet)}, \aa^\bullet) = \rhat(\aa^{(\bullet)}, \overline{\aa^\bullet}) = \rho(\aa^{(\bullet)}, \overline{\aa^\bullet})$ by \cite{DFMS19}.

\begin{corollary} \label{cor.rhatV}
    Let $\aa$ and $\bb$ be squarefree monomial ideals in $R=\kk[x_1,\ldots,x_n].$ Then, $$\rhat(\aa^{(\bullet)}, \overline{\bb^\bullet})  = \dfrac{1}{\min\{ \langle \mathbf{a}, \mathbf{b} \rangle \mid \mathbf{a}\in \SP(\aa) ,  \mathbf{b}\in \SP(\bb^\vee)\}} = \rhat((\bb^\vee)^{(\bullet)}, \overline{(\aa^\vee)^\bullet}). $$ In particular, for any squarefree monomial ideal $\aa \subseteq R$, $$\rhat(\aa^{(\bullet)}, \overline{\aa^\bullet})  = \rhat((\aa^\vee)^{(\bullet)}, \overline{(\aa^\vee)^\bullet}). $$
\end{corollary}

\begin{proof}
    The proof follows from Corollary \ref{cor:asym-res-IP}, the facts that $\Delta(\aa^{(\bullet)}) = \SP(\aa)$, $(\aa^\vee)^\vee=\aa$, and that both $\SP(\aa), \SP(\bb^\vee)$ are polyhedra.
\end{proof}

%%%%%%% generalization %%%%%%%%%%

\section{Invariant ideals, Rees packages and convex bodies}
\label{sec.Rees}

In this section, we attempt to understand combinatorially the asymptotic resurgence from associated convex bodies, and look for statements that are similar to Theorem \ref{thm.rhohatViaNO} for more general classes of ideals. Throughout the section, we will assume that $R$ is a Nagata algebra over a field $\kk$. Our approach is based on the notion of \emph{Rees package} introduced and studied recently in \cite{BDHM}.

%%%%%%%%% Rees packages %%%%%%%%%%

\begin{definition}[\protect{\cite[Definition 3.2]{BDHM}}]
	Fix a $\kk$-vector space basis $\B = \{\b_i ~\big|~ i \in \Lambda\}$ for $R$. We say that $v : R \rightarrow \ZZ_{\ge 0}$ is a \emph{$\B$-monomial function} if for every $f = \sum_{i \in \Lambda} c_i\b_i \in R$, we have $v(f) = \min\{v(\b_i) ~\big|~ c_i \not= 0\}.$ A valuation of $R$, that is also a $\B$-monomial function, is called a \emph{$B$-monomial valuation}.
\end{definition}
Recall that a hyperplane in $\RR^d$ defined by $\langle \h, \x\rangle = c$, where $\h = (h_1, \dots, h_d) \in \RR^d$, $\x = (x_1, \dots, x_d)$ and $c \in \RR$ is called \emph{non-coordinate} if $c \not= 0$.

\begin{definition}[\protect{\cite[Proposition 3.4 and Definition 3.5]{BDHM}}]
	\label{def.Rees}
	Let $v_1, \dots, v_d$ be $\B$-monomial functions of $R$ and set $\v = (v_1, \dots, v_d)$. Let $\Gamma \subseteq \RR^d_{\ge 0}$ be an integral polyhedron such that for every $s \in \NN$, the integral closure $\overline{I^s}$ is spanned as a $\kk$-vector space by the set
	$$\{\b \in \B ~\big|~ \v(\b) \in s\Gamma\}.$$
	Let $H_1, \dots, H_r$ be the non-coordinate supporting hyperplanes of $\Gamma$, where $H_k$ is given by $\langle \h_k, \x\rangle = c_k $ for some $\h_k \in \ZZ^d_{\ge 0}$ and $c_k\in \ZZ_{> 0}$. We call $(\B,\v,\Gamma)$ a \emph{Rees package} for $I$ if the $\B$-monomial function $V_k : = \langle \h_k, \v\rangle$ is a valuation of $R$ for any $k = 1, \dots, r$. In this case, $\{V_1, \dots, V_r\}$ are precisely the \emph{Rees valuations} of $I$ and for any $u \in \QQ_{\ge 0}$, the \emph{rational power} $\overline{I^u}$ is spanned as a $\kk$-vector space by
	$$\{\b \in \B ~\big|~ \v(\b) \in u\Gamma\}.$$
	See, for example, \cite{SH} for the definitions of Rees valuations and rational powers. Rees packages were shown to exist for several classes of monomial and invariant ideals in \cite{BDHM}.
\end{definition}

Our next definition constructs a closed set associated to a graded family of ideals with respect to an ideal with a given Rees package.

\begin{definition}
	\label{def.convexBD}
	Let $\bb \subseteq R$ be an ideal and assume that $\R = (\B,\v,\Gamma)$ is a Rees package for $\bb$. For any graded family $\aa_\bullet = \{\aa_i\}_{i \ge 1}$, define
	$$V_\R(\aa_i) = \conv \left( \{\v(\b) ~\big|~ f \in \aa_i \text{ and } \b \text{ is a } \text{$\B$-monomial in } f\} \right) \subseteq \RR^d_{\ge 0},$$
 and set
	$$\Gamma_\R(\aa_\bullet) = \overline{\bigcup_{k = 1}^\infty \dfrac{1}{k}V_\R(\aa_k)}.$$
\end{definition}

\begin{remark} \label{rmk.convex} 
	We say that $\v$ is a \emph{vector of valuations} if the $\B$-monomial functions $v_1, \dots, v_d$ are valuations of $R$. In this case, the condition that $V_k$ is a valuation, for $k = 1, \dots, r$, for $(\B,\v,\Gamma)$ to be a Rees package in Definition \ref{def.Rees} is automatically satisfied. It is not hard to see that if $\v$ is a vector of valuations then $\Gamma_R(\aa_\bullet)$ is a closed convex set. This is the case in all known examples where Rees packages exist; see \cite{BDHM}.
\end{remark} 

We shall use an example given in \cite{BDHM} to illustrate the construction introduced in Definition \ref{def.convexBD}.

\begin{example} \label{ex:sym-matrix}
Let $Y=[y_{i,j}]$ be a $m \times m$ symmetric matrix consists of indeterminates $y_{i,j}$ with $1 \le i,j \le m$ and $m \ge 3$. Let $R=\kk[y_{i,j} ~\big|~ i,j \in [m]]$ and let $\bb=I_{m-1}(Y)$ be the ideal generated by all $(m-1) \times (m-1)$-minors of $Y$ in $R$. Let $\B$ be the set of standard monomials in $R$, see \cite[Setup 3.16]{BDHM}.

It follows from \cite[Theorem 3.17]{BDHM} that $\bb$ admits a Rees packages $\R = (\B,\v,\Gamma)$, where $\v = (\gamma_1, \dots, \gamma_m)$ is a vector of $\B$-monomial valuations, and
$$\Gamma =\text{convexhull} \{\left(m-1,m-2,\ldots, 1, 0\right)\} +\RR^m_{\ge 0}.$$

Consider the graded family $\aa_\bullet= \bb^{(\bullet)}$.  By \cite[Proposition 4.3]{JeffriesMontanoVarbaro15}, for all $s \in \NN$, $\bb^{(s)}$ is generated by standard monomials,  $\bb^{(2k)}=\left( \bb^{(2)}\right)^k$ and $\bb^{(2k+1)}=\bb\left( \bb^{(2)}\right)^k$ for all $k$. Also, $\bb^{(2)}=\bb^2+I_m(Y).$ Thus, for all $k \ge 1$,
$$\v \left(\bb^{(2k)}\right) =k\v\left(\bb^{(2)}\right) \text{ and } \v \left(\bb^{(2k+1)}\right) =k\v\left(\bb^{(2)}\right)+\v(\bb).$$

By the construction of $\v$ given in \cite[Theorem 3.17]{BDHM}, we get
$$\v(\bb)= (m-1,m-2,\ldots, 1, 0) \text{ and } \v(I_m(Y))= (m,m-1,\ldots, 2, 1).$$
Therefore,
$$V_\R(\aa_{2k}) = \conv \left( \{ k(m,m-1,\ldots,2,1), 2k(m-1,m-2,\ldots,1,0)\} \right) +\RR_{\ge 0}^m,$$
and
\begin{align*}
	V_{\R}(\aa_{2k+1}) =  \conv & \left\{ (k+1)(m,m-1,\ldots,2,1)-(1,1,\ldots,1) , \right. \\
	&  \quad \left.  (2k+1)(m-1,m-2,\ldots,1,0) \right\} +\RR_{\ge 0}^m.
\end{align*}
	
Thus, taking the limit gives us
	\begin{align*} \Gamma_\R(\aa_\bullet) & = \overline{\bigcup_{k = 1}^\infty \dfrac{1}{k}V_\R(\aa_k)} \\& =\conv \left( \left\{ \frac{1}{2}(m,m-1,\ldots,2,1), (m-1,m-2,\ldots,1,0)\right\} \right) +\RR_{\ge 0}^m .\end{align*}
\end{example}

%%%%%%%%%% generalization for one ideal %%%%%%%%%%%

Our next result establishes $\rhat(\aa_\bullet, \overline{\bb_\bullet})$ when $\bb_\bullet = \bb^{\bullet}$ is the family of ordinary powers of an ideal $\bb \subseteq R$ with Rees packages.

\begin{theorem}
	\label{thm.rhat_b}
	Let $\bb \subseteq R$ be a nonzero ideal with Rees package $\R = (\B,\v,\Gamma)$, in which $\v$ is a vector of valuations. Let $\aa_\bullet = \{\aa_i\}_{i \ge 1}$ be a graded family of nonzero ideals in $R$. Then,
	$$\rhat(\aa_\bullet, \overline{\bb^\bullet}) = \sup \{\lambda > 0 \mid \lambda \cdot \Gamma_\R(\aa_\bullet) \not\subseteq \Gamma\}.$$
\end{theorem}

\begin{proof}
	Consider $s,r \in \NN$ such that $\aa_{st} \not\subseteq \overline{\bb^{rt}}$ for $t \gg 0$. It follows that, for $t \gg 0$, $\aa_{st}$ is not spanned by the $\B$-monomials whose $\v$-values are in $rt\Gamma$. Thus, for $t \gg 0$, there exists $f \in \aa_{st}$ and a $\B$-monomial $\b$ in $f$ such that $\v(\b) \not\in rt\Gamma$. Since $\Gamma$ is closed and convex, this implies that
	$$\frac{1}{st}V_{\R}(\aa_{st}) \not\subseteq \frac{rt}{st}\Gamma \text{ for } t \gg 0.$$
	That is, $\frac{s}{r} \cdot \frac{1}{st}V_{\R}(\aa_{st}) \not\subseteq \Gamma$ for $t \gg 0$. Therefore, $$\frac{s}{r} \cdot \Gamma_\R(\aa_\bullet) \not\subseteq \Gamma.$$
	
	Conversely, let $s,r \in \NN$ be such that $\frac{s}{r} \Gamma_\R(\aa_\bullet) \not\subseteq \Gamma$. Since $\Gamma$ is closed and convex, it follows that there exists $k \in \NN$ with
	$$\frac{s}{r}\cdot \frac{1}{k}V_{\R}(\aa_k) \not\subseteq \Gamma,$$
	i.e., $s\cdot V_{\R}(\aa_k) \not\subseteq rk\Gamma$. Since $\v$ is a vector of valuations and $\aa_k^{st} \subseteq \aa_{kst}$, we have $stV_{\R}(\aa_k) \subseteq V_{\R}(\aa_{kst})$ for all $t \in \NN$. This implies that, for all $t \in \NN$, $V_{\R}(\aa_{kst}) \not\subseteq rkt\Gamma$, and so there exists $f \in \aa_{kst}$ and $\b$ a $\B$-monomial in $f$ such that $\v(\b) \not\in rkt\Gamma$. Particularly, $\aa_{kst} \not\subseteq \overline{\bb^{rkt}}$ for all $t \in \NN$. Hence,
	$$\rhat(\aa_\bullet, \overline{\bb^\bullet}) \le \frac{s}{r}.$$
	The result is established.
\end{proof}

%%%%%%%% generalization for Noetherian b-family %%%%%%%%%%%

Theorem \ref{thm.rhat_b} has a generalization when $\bb_\bullet$ is a Noetherian family. For a graded family $\bb_\bullet = \{\bb_i\}_{i \ge 1}$ and a positive integer $k$, set
$$\R^{[k]}(\bb_\bullet) = \bigoplus_{i \ge 0} \bb_{ik} t^i \subseteq R[t].$$
This is the $k$-th Veronese subalgebra of the Rees algebra $\R(\bb_\bullet)$.

\begin{theorem}
	\label{thm.rhat_Noetherian}
	Let $\aa_\bullet$ and $\bb_\bullet$ be graded families of nonzero ideals in $R$. Assume that for some $k \in \NN$, $\R^{[k]}(\bb_\bullet)$ is a standard graded $R$-algebra, and let $\bb_k^\bullet$ denote the graded family $\{\bb_k^i\}_{i \ge 1}$. Suppose that $\bb_k$ has a Rees package $\R_k = (\B_k,\v_k,\Gamma_k)$, in which $\v_k$ is a vector of valuations. Then,
	$$\rhat(\aa_\bullet, \overline{\bb_\bullet}) = \dfrac{1}{k} \sup \{ \lambda > 0 \mid \lambda \cdot \Gamma_{\R_k}(\aa_\bullet) \not\subseteq \Gamma_k\}.$$
\end{theorem}

\begin{proof}
	By \cite[Lemma 4.14]{HKNN23}, we have
	$$\rhat(\aa_\bullet, \overline{\bb_\bullet}) = \dfrac{1}{k} \rhat(\aa_\bullet, \overline{\bb_k^\bullet}).$$
	The result then follows from Theorem \ref{thm.rhat_b}.
\end{proof}

%%%%%%%%%%% generalization for b-equivalent %%%%%%%%%%%%%

Recall that a graded family $\bb_\bullet = \{\bb_i\}_{i \ge 1}$ is said to be $\bb$-equivalent for an ideal $\bb \subseteq R$ if there exists a positive integer $k$ such that for any $i \ge 1$,
$$\bb_{i+k} \subseteq \bb^i \subseteq \bb_i.$$
Our next result computes $\rhat(\aa_\bullet, \bb_\bullet) = \rhat(\aa_\bullet, \overline{\bb_\bullet})$ when $\bb_\bullet$ is $\bb$-equivalent.

\begin{theorem}
	\label{thm.rhat_b-equivalent}
	Let $\bb \subseteq R$ be a nonzero ideal with Rees package $\R = (\B,\v,\Gamma)$, in which $\v$ is a vector of valuations. Let $\aa_\bullet = \{\aa_i\}_{i \ge 1}$ be a filtration and let $\bb_\bullet = \{\bb_i\}_{i \ge 1}$ be a graded family of nonzero ideals in $R$. Suppose that $\bb_\bullet$ is $\bb$-equivalent. Then,
	$$\rhat(\aa_\bullet, \bb_\bullet) = \rho(\aa_\bullet, \overline{\bb_\bullet}) = \rhat(\aa_\bullet, \overline{\bb_\bullet}) = \sup \{\lambda > 0 : \lambda \cdot \Gamma_\R(\aa_\bullet) \not\subseteq \Gamma\}.$$
\end{theorem}

\begin{proof}
	The first two equalities follow from \cite[Theorem 4.8]{HKNN23}. Also, as observed in \cite[Theorem 4.8]{HKNN23}, since $\bb_\bullet$ is $\bb$-equivalent, by \cite[Corollary 3.3]{HKNN23}, we have
	$$\rhat(\aa_\bullet, \overline{\bb_\bullet}) = \rhat(\aa_\bullet, \overline{\bb^\bullet}).$$
	The result now follows from Theorem \ref{thm.rhat_b}.
\end{proof}

Here, we illustrate  Theorems \ref{thm.rhat_b} and \ref{thm.rhat_b-equivalent} using Example \ref{ex:sym-matrix}.
\begin{example}
	As we saw in Example \ref{ex:sym-matrix}, $\Gamma =\text{convexhull} \{\left(m-1,m-2,\ldots, 1, 0\right)\} +\RR^m_{\ge 0}$  and
	\begin{align*}
		\Gamma_\R(\aa_\bullet) & = \overline{\bigcup_{k = 1}^\infty \dfrac{1}{k}V_\R(\aa_k)} \\
		& =\conv \left( \left\{ \frac{1}{2}(m,m-1,\ldots,2,1), (m-1,m-2,\ldots,1,0)\right\} \right) +\RR_{\ge 0}^m .
	\end{align*}

Observe that $\Gamma \subseteq \Gamma_\R(\aa_\bullet)$ and the supremum value of $\lambda$ such that $\lambda \cdot \Gamma_\R(\aa_\bullet) \not\subseteq \Gamma$ is determined by where the line connecting the origin and $\frac{1}{2}\left(m, m-1, \dots, 2,1\right)$ first intersects $\Gamma$. This point is computed to be $\left(m-1, \frac{(m-1)^2}{m}, \frac{(m-1)(m-2)}{m}, \dots, \frac{m-1}{m}\right)$. Therefore,
$$\sup \left\{\lambda > 0 ~\big|~ \lambda \cdot \Gamma_\R(\aa_\bullet) \not\subseteq \Gamma\right\} = \frac{2(m-1)}{m}.$$
 Hence, due to Theorem \ref{thm.rhat_b}, we get
	$$\rhat(\aa_\bullet, \overline{\bb_\bullet}) = \frac{2(m-1)}{m}.$$
Note that in this case, it is actually the usual asymptotic resurgence of the ideal $\bb$, which has been recently computed in \cite[Section 4]{KM24}.
\end{example}

Coupling Theorems \ref{thm.rhat_b}, \ref{thm.rhat_b-equivalent}, and the work in \cite{BDHM} on Rees packages, we obtain the following corollary.

\begin{corollary}
	\label{cor.ReesPackage}
	Let $\bb \subseteq R$ be an ideal that belongs to the following classes:
	\begin{enumerate}
		\item[(i)] monomial ideals in affine semigroup rings;
		\item[(ii)] sums of products of determinantal ideals of generic matrices;
		\item[(iii)] sums of products of determinantal ideals of symmetric matrices;
		\item[(iv)] sums of products of ideals of Pfaffians of skew-symmetric matrices;
		\item[(v)] products of determinantal ideals of Hankel matrices.
	\end{enumerate}
	Let $\R = (\B,\v,\Gamma)$ be a Rees package of $\bb$. Let $\aa_\bullet$ and $\bb_\bullet$ be graded families of ideals in $R$. Assume either
	\begin{enumerate}
		\item $\bb_\bullet = \bb^\bullet$ is the family of ordinary powers of $\bb$; or
		\item $\aa_\bullet$ is a filtration and $\bb_\bullet$ is $\bb$-equivalent.
	\end{enumerate}
	Then,
	$$\rhat(\aa_\bullet, \overline{\bb_\bullet}) = \sup \{\lambda > 0 ~\big|~ \lambda \cdot \Gamma_\R(\aa_\bullet) \not\subseteq \Gamma\}.$$
\end{corollary}

%%%%%%%%%%%%%%%%%%%%%%%%%%%%%%
%%%%%%%%%%%%%%%%%%%%%%%%%%%%%%

\section{Truncation of graded families} \label{sec.Truncate}

In this section, we investigate whether resurgences of not necessarily Noetherian graded families can be approximated by those of Noetherian families. Particularly, we will look at resurgences of \emph{truncated} families. We shall begin with a necessary definition.

\begin{definition} \label{def.truncate}
	Let $\aa_\bullet$ be a graded family of ideals in $R$. The \emph{$n$-th truncation of} $\aa_\bullet$, denoted by $\aa_{n,\bullet} = \{\aa_{n,k}\}_{k \in \NN}$, is a graded family of ideals in $R$, where \[   \aa_{n,k}=  \left\{
	\begin{array}{ll}
		\aa_k & \text{ if } k \leq n \\
		\sum\limits_{i,j > 0,\; i+j=k} \aa_{n,i} \aa_{n,j} & \text{ if } k > n. \\
	\end{array}
	\right. \]
\end{definition}

Our next result shows that truncation works well for understanding the so-called skew Waldschmidt constant of a graded family of ideals.  Observe that if $\aa_\bullet$ is a graded family of ideals in a domain $R$ and $v$ is a valuation of the fraction field of $R$, then $\{v(\aa_n)\}_{n \in \NN}$ is subadditive, and so the following limit exists and is called the \emph{skew Waldschmidt} constant of $\aa_\bullet$ with respect to $v$:
$$\vhat(\aa_\bullet) = \lim\limits_{n \rightarrow \infty} \frac{v(\aa_n)}{n} = \inf_{n \in \NN} \frac{v(\aa_n)}{n}.$$
The valuation $v$ is said to be supported on $R$ if $v(x) \geq 0$ for all $x \in R$.

\begin{theorem}
\label{thm.vhatTruncation}
Let $\aa_\bullet$  be a graded family of ideals in a domain $R$.
\begin{enumerate}
	\item If $v$ is a valuation of fraction field of $R$ that is supported on $R$. Then,
    $$\vhat(\aa_\bullet) =\inf\limits_{n \ge 1 } \vhat(\aa_{n,\bullet})=\lim_{n \to \infty} \vhat(\aa_{n,\bullet}).$$
    \item  If $R$ is a polynomial ring and $\aa_\bullet$  is a graded family of monomial ideals, then
    $$\Delta(\aa_\bullet) =\overline{\bigcup_{n =1}^{\infty} \Delta(\aa_{n,\bullet})}.$$
\end{enumerate}
\end{theorem}

\begin{proof}
Consider any integers $n \ge m \ge 1$. Note that $\aa_{m,k} \subseteq \aa_{n,k} \subseteq  \aa_k$, for every $k \ge 1$. Therefore, $v(\aa_k) \le v(\aa_{n,k}) \le v(\aa_{m,k})$  for every $k \ge 1$. This implies that $$\vhat(\aa_\bullet) \le \vhat(\aa_{n,\bullet}) \le \vhat(\aa_{m,\bullet}) \text{ for all } n\ge m \ge 1.$$
Furthermore, $\{\vhat(\aa_{n,\bullet})\}_{n\ge 1}$ is a non-increasing sequence, hence,
$$\vhat(\aa_\bullet) \le \inf\limits_{n \ge 1 } \vhat(\aa_{n,\bullet})=\lim\limits_{n \to \infty} \vhat(\aa_{n,\bullet}).$$

On the other hand, let $\epsilon >0$ be any positive real number.  Then, there exists $k_0$ such that, for all $k \ge k_0$,
$$\frac{v(\aa_k)}{k}< \vhat(\aa_\bullet) +\epsilon.$$
For $n\ge k_0$, we then have
$$\vhat(\aa_{n,\bullet}) \le \frac{v(\aa_{n,n})}{n} =\frac{v(\aa_n)}{n} <\vhat(\aa_\bullet) +\epsilon.$$
This is true for an $\epsilon > 0$, so $\lim\limits_{n \to \infty} \vhat(\aa_{n,\bullet}) \le \vhat(\aa_\bullet).$ We have established the equality
$$\lim\limits_{n \to \infty} \vhat(\aa_{n,\bullet}) = \vhat(\aa_\bullet).$$

To prove the second statement, again, consider any integers $n \ge m \ge 1$. Since $\aa_{m,k} \subseteq \aa_{n,k} \subseteq  \aa_k$, for every $k \ge 1$,  we have
$$\NP(\aa_{m,k}) \subseteq \NP(\aa_{n,k}) \subseteq \NP( \aa_k).$$
Consequently, $$\Delta(\aa_{m,\bullet}) \subseteq \Delta(\aa_{n,\bullet}) \subseteq \Delta(\aa_\bullet).$$ Thus, $\Delta(\aa_{1,\bullet}) \subseteq \Delta(\aa_{2,\bullet}) \subseteq \cdots $ is an ascending sequence of closed sets.

It is easy to note that
$$\overline{\bigcup_{n =1}^{\infty} \Delta(\aa_{n,\bullet})} \subseteq \Delta(\aa_\bullet) .$$
Now consider any $\lambda \in \Delta(\aa_\bullet)$. Then, there exists a sequence $ \left\{ \frac{\alpha_{n_k}}{n_k} \right\}_{k \ge 1}$ with $\alpha_{n_k} \in \NP(\aa_{n_k})$ such that  $\lim\limits_{k \to \infty} \frac{\alpha_{n_k}}{n_k} = \lambda.$ Note that for every $k \ge 1$,
$$\frac{\alpha_{n_k}}{n_k}  \in \frac{1}{n_k}\NP(\aa_{n_k, n_k}) \subseteq \Delta(\aa_{n_k,\bullet}).$$
Thus,
$$\frac{\alpha_{n_k}}{n_k} \in \overline{\bigcup_{n =1}^{\infty} \Delta(\aa_{n,\bullet})} \text{ for all } k \ge 1.$$
Since $\overline{\bigcup_{n =1}^{\infty} \Delta(\aa_{n,\bullet})} $ is a closed set, we get that $\lambda \in \overline{\bigcup_{n =1}^{\infty} \Delta(\aa_{n,\bullet})}.$  Hence, the assertion follows.
\end{proof}

Our last main result of the paper is stated as follows.

\begin{theorem}
    \label{thm.resurgencetruncation1}
    Let $\aa_\bullet$ and $\bb_\bullet$ be graded families of ideals in a Noetherian commutative ring $R$. For $n \in \NN$, let $\aa_{n,\bullet}$ denote the $n$-th truncation graded family of $\aa_\bullet$. We have
    $$\rho(\aa_\bullet,\bb_\bullet) =\sup\limits_{n \ge 1 } \rho(\aa_{n,\bullet},\bb_\bullet)=\lim_{n \to \infty} \rho(\aa_{n,\bullet},\bb_\bullet) \text{ and }
    \rhat(\aa_\bullet,\bb_\bullet) =\sup\limits_{n \ge 1 } \rhat(\aa_{n,\bullet},\bb_\bullet)=\lim_{n \to \infty} \rhat(\aa_{n,\bullet},\bb_\bullet).$$
\end{theorem}

\begin{proof}
    Since $\aa_{n,\bullet} \subseteq \aa_\bullet$, we have $\rho(\aa_{n,\bullet},\bb_\bullet) \le \rho(\aa_\bullet,\bb_\bullet)$ for all $n\in \NN$. Moreover, since $\aa_{n,\bullet} \subseteq \aa_{n+1,\bullet}$, the sequence $\{\rho(\aa_{n,\bullet},\bb_\bullet)\}_{n\in \NN}$ is a non-decreasing sequence. Thus, the limit $\lim_{n \to \infty} \rho(\aa_{n,\bullet},\bb_\bullet)$ exists and
    $$\rho(\aa_\bullet,\bb_\bullet) \ge \sup\limits_{n \ge 1 } \rho(\aa_{n,\bullet},\bb_\bullet)=\lim_{n \to \infty} \rho(\aa_{n,\bullet},\bb_\bullet).$$

    Suppose that $\rho(\aa_\bullet,\bb_\bullet) > \sup\limits_{n \ge 1 } \rho(\aa_{n,\bullet},\bb_\bullet)$, then there exists $s,r \in \NN$ such that
    $$\rho(\aa_\bullet,\bb_\bullet) \ge \frac{s}{r} > \sup\limits_{n \ge 1 } \rho(\aa_{n,\bullet},\bb_\bullet) \text{ and } \aa_{s} \not\subseteq \bb_{r}.$$
    From $\frac{s}{r} > \sup\limits_{n \ge 1 } \rho(\aa_{n,\bullet},\bb_\bullet)$, we have $\aa_{n,s}\subseteq \bb_r$, for every $n\in \NN$. In particular, $\aa_s = \aa_{s,s}\subseteq \bb_r$, which is a contradiction. Therefore, $\rho(\aa_\bullet,\bb_\bullet) = \sup\limits_{n \ge 1 } \rho(\aa_{n,\bullet},\bb_\bullet)$. The statement for $\rhat$ is proved identically.
\end{proof}

We end the paper with an example exhibiting that a similar statement to that of Theorem \ref{thm.resurgencetruncation1} may not hold if truncations of $\bb_\bullet$ are taken instead of those of $\aa_\bullet$.

\begin{example}
\label{ex.rhoTruncation}
    In general, $\displaystyle \lim_{n \to \infty} \rho(\aa_\bullet,\bb_{n,\bullet})$ and $\rho(\aa_\bullet,\bb_\bullet)$ may be different (and similarly, $\displaystyle \lim_{n \to \infty} \rhat(\aa_\bullet,\bb_{n,\bullet})$ and $\rhat(\aa_\bullet,\bb_\bullet)$ may be different, etc.). Consider the families $\aa_\bullet$ and $\bb_\bullet $ with $\aa_i = I^i$ and
		$$\bb_i = \left\{ \begin{array}{lll} I & \text{if} & i \not= 2 \\ I^2 & \text{if} & i = 2.\end{array}\right.$$
		
By \cite[Example 2.10 (2)]{HKNN23}, we have
		$$\rho(\aa_\bullet,\bb_\bullet)=\dfrac{1}{2} \text{ and } \rhat(\aa_\bullet,\bb_\bullet) = \rhat(\aa_\bullet, \overline{\bb_\bullet}) =-\infty < \dfrac{\vhat(\bb_\bullet)}{\vhat(\aa_\bullet)}.$$

 For $n\ge 5$, by induction on $k$, we have the following formula for the elements of the $n$-th truncation $\bb_{n,\bullet}$:
    $$\bb_{n,k} = \left\{ \begin{array}{lll} I^{\lceil \frac{k}{n} \rceil} & \text{if} & k \not= 2 \\ I^2 & \text{if} & k = 2.\end{array}\right.$$
    Hence, for each $n \ge 5$,
    \[
    \rho(\aa_\bullet,\bb_{n,\bullet}) = \sup \left \{ \frac{s}{r} ~ \Big| ~ I^s \not \subseteq I^{\lceil \frac{r}{n} \rceil} \right \} = \sup_s \left \{ \frac{s}{sn+1} \right \} = \frac{1}{n}.
    \]
    Therefore, $\displaystyle \lim_{n \to \infty} \rho(\aa_\bullet,\bb_{n,\bullet}) = 0 \neq \rho(\aa_\bullet,\bb_\bullet)=\dfrac{1}{2}$.

    On the other hand, for each $n \ge 5$,
    \[
    \rhat(\aa_\bullet,\bb_{n,\bullet}) = \sup \left \{ \frac{s}{r} ~ \Big| ~ I^{st} \not \subseteq I^{\lceil \frac{rt}{n} \rceil}, t\gg 1 \right \} = \sup_s \left \{ \frac{s}{sn+1} \right \} = \frac{1}{n}.
    \]
    The second equality holds as $I^{st} \not \subseteq I^{\lceil \frac{rt}{n} \rceil}$ implies that $\lceil \frac{rt}{n} \rceil \ge st+1$, hence, $rt\ge stn+t$. Therefore, $\displaystyle \lim_{n \to \infty} \rhat(\aa_\bullet,\bb_{n,\bullet}) = 0 \neq \rhat(\aa_\bullet,\bb_\bullet)=-\infty$.

    Furthermore, if we pick $I$ to be a normal ideal, then $\overline{\bb_{n,i}}=\bb_{n,i}$ and $\overline{\bb_{i}}=\bb_{i}$, so
    $$\displaystyle \lim_{n \to \infty} \rhat(\aa_\bullet,\overline{\bb_{n,\bullet}}) = \lim_{n \to \infty} \rhat(\aa_\bullet,\bb_{n,\bullet}) \neq \rhat(\aa_\bullet,\bb_\bullet) = \rhat(\aa_\bullet,\overline{\bb_\bullet}).$$
Note that in this case if $I$ is a monomial ideal, then we still have
$$\rhat(\aa_\bullet, \overline{\bb_\bullet}) =\sup \{ \lambda ~\big|~ \lambda \cdot \Delta(\aa_\bullet) \not\subseteq \Delta(\bb_\bullet), \lambda >0 \} .$$
In fact, we see that $\Delta(\aa_\bullet) = \NP(I)$ and that $$\Delta(\bb_\bullet) = \overline{\bigcup_{k=1}^\infty \frac{1}{k}\NP(\bb_k)} = \overline{\bigcup_{k=1}^\infty \frac{1}{k}\NP(I)}=\RR^n_{\ge 0}.$$
 Thus,
 $$\sup \{ \lambda > 0 ~\big|~ \lambda \cdot \Delta(\aa_\bullet) \not\subseteq \Delta(\bb_\bullet), \lambda >0 \} = -\infty = \rhat(\aa_\bullet, \overline{\bb_\bullet}).$$
\end{example}

%%%%%%%%%%%%%%%%%%%%%%%%%%%%%%%%%%%%%%%%%%

\noindent\textbf{Acknowledgement.}
Part of this work was done while the second and third authors visited Tulane University. We thank Tulane University for its supports and hospitality. The first author is partially supported by Simons Foundation. The second author acknowledges the partial support by the Science and Engineering Research Board (SERB) MATRICS Grant (Grant Number: MTR/2023/000335). The third author acknowledges the support by AMS-Simons Travel Grant. The last author is supported by the Vietnam National Foundation for Science and Technology Development (NAFOSTED).

\vskip 2mm\noindent{\bf Disclosure statement:}
The authors have no competing interests to declare that are relevant to the content of this article.
%%%%%%%%%%%%%%%%%%%%%%%%%%%%%%%%
%%%%%%%%%%%%%%%%%%%%%%%%%%%



\end{document}